\documentclass[11pt,a4]{article}
\usepackage{amssymb,amsmath,amsthm,enumerate}
\usepackage{psfrag}
\usepackage{dcolumn}
\usepackage{color}
\usepackage{bm}
\usepackage{color}

\begin{document} 

\title{ The Stationary Boltzmann equation for a two component gas in the
  slab with different molecular masses.\\}
\author{
St\'ephane Brull *\\
\hspace{1.in}\\
\hspace{1.in}\\
\hspace{1.in}\\}
\date{}
\maketitle
\newtheorem{remark}{Remark}
\newtheorem{theorem}{Theorem}
\newtheorem{proposition}{Proposition}
\newtheorem{lemma}{Lemma}
\newtheorem{property}{Property}
\newtheorem{example}{Example}
\newtheorem{definition}{Definition}
\newtheorem{corollary}[theorem]{Corollary}
\renewcommand{\thetheorem}{\thesection.\arabic{theorem}}
\renewcommand{\thelemma}{\thesection.\arabic{lemma}}
\renewcommand{\theproperty}{\thesection.\arabic{property}}
\renewcommand{\theexample}{\thesection.\arabic{example}}
\renewcommand{\thedefinition}{\thesection.\arabic{definition}}
\renewcommand{\thecorollary}{\thesection.\arabic{corollary}}
\renewcommand{\theequation}{\thesection.\arabic{equation}}
\newcommand{\sectionnew}[1]{
\section{#1}\setcounter{equation}{0}
\setcounter{theorem}{0}
\setcounter{lemma}{0}}
\begin{abstract}

The stationary Boltzmann equation for hard and soft forces in the context
of a two component gas is considered in the slab when the 
molecular masses of the 2 component are different. An $L^{1}$
existence theorem is proved when one component satisfies a given
indata profile and the other component satisfies diffuse
reflection at the boundaries. Weak $L^{1}$ compactness is
extracted from the control of the entropy production term. 
%{\it To cite this article: A. Nom1, A. Nom2, C. R.
%Acad. Sci. Paris, Ser. I 340 (2005).}
\end{abstract}
%\end{frontmatter}
\begin{tabular}{p{13.7cm}}
\hline
{\small *Mathematiques appliquees de Bordeaux, University of Bordeaux
  I, 
351 cours de la Lib\'eration 33405 Talence Cedex, France.}
\end{tabular}
% Maintenant la version abr�g�e en anglais, si pr�sente
%\selectlanguage{english}
\setcounter{equation}{0}
\section{Introduction and setting of the problem.}
This article is devoted to the proof of an existence theorem for 
the stationary Boltzmann equation in the situation of a two component 
gas having  different molecular masses for the geometry 
of the slab. The slab being represented by the interval $[-1,1]$, the Boltzmann equation reads 
\begin{eqnarray} \label{cA} 
\xi \frac{\partial}{\partial x} f_A(x,v) = 
Q_{AA}(f_A,f_A)(x,v) + Q_{AB}(f_A,f_B)(x,v) ,
\\ \label{cB}\xi \frac{\partial}{\partial x} f_B(x,v) = 
Q_{BB}(f_B,f_B)(x,v) + Q_{BA}(f_B,f_A)(x,v) ,
\\ \nonumber x \in [-1,1], \hspace*{1 mm} v \in \mathbb{R}^3.
\end{eqnarray}
The non-negative functions represent the distribution functions $f_A$ 
and $f_B$ of the $A$ and the $B$ component and $\xi$ is the velocity
component in the $x$ direction.
For for any $\alpha , \,  \beta \, \in \{ A,B\}$, $Q_{\alpha,
  \beta}$ corresponds to the Boltzmann collision operator 
between the species $\alpha$ and $\beta$. It is defined by 
\begin{eqnarray} \label{col}
Q_{\alpha, \beta}(v) = \int_{\mathbb{R}^3 \times \mathcal{S}^2} 
\mathcal{B}^{\alpha , \beta} \left(  f_{\alpha}(x,v_{\ast}^{\prime}) 
 f_{\beta}(x,v^{\prime}) -   f_{\beta}(x,v_{\ast}) 
 f_{\alpha}(x,v) \right) d \omega dv_{\ast} \hspace*{1 mm}
\end{eqnarray}
with
\begin{eqnarray} \label{velo}
\nonumber v^{\prime(\beta \alpha)} = v + \frac{2 m^{\beta}}{m^{\alpha}+m^{\beta}}
\langle v_{\ast} -v  , \omega  \rangle \omega
, \hspace*{5 mm} v_{\ast}^{\prime(\beta \alpha)} = v_{\ast} - \frac{2 m^{\beta}}{m^{\alpha}+m^{\beta}}
\langle v_{\ast} -v  , \omega  \rangle \omega .
\\ 
\end{eqnarray}
In the formula (\ref{velo}), $ v^{\prime(\beta \alpha)}$ and  
$v_{\ast}^{\prime(\beta \alpha)}$ represent the post-colisional 
velocities between the species $\alpha$ and $\beta$ and $m^{\alpha}$
represents the mass of the specy $\alpha$.
For more precisions on the model we refer to (\cite{chap}, \cite{ABT}).

 $ \langle \cdot \, \, ,\cdot \, \rangle$ denotes the Euclidean inner product in
$\mathbb{R}^{3}$. Let $\omega$ be represented by the polar angle
(with polar axis along $v - v_{\ast}$) and the azimutal angle
$\phi$.

For the sake of clarity, recall the invariant properties 
of the collision operator $Q_{\alpha,\beta}$, $\{\alpha ,\beta\} \in
\{A,B\}$. For more details we refer to (\cite{Do}).

\begin{property} \label{cons}
For $\alpha , \beta$ $\in $ $\{A,B\}$, with $\alpha \not= \beta $, it
holds that
\begin{eqnarray*}
\int_{\mathbb{R^3}} (1,m^{\alpha}v , m^{\alpha}|v|^2) Q_{\alpha ,
  \alpha}
(f_{\alpha}, f_{\alpha}) dv =0,
\\ \int_{\mathbb{R^3}}  Q_{\alpha ,
  \beta}
(f_{\alpha}, f_{\beta}) dv =0,
\\ \int_{\mathbb{R^3}}  m^{\alpha} v \, Q_{\alpha ,
  \beta}
(f_{\alpha}, f_{\beta}) dv + \int_{\mathbb{R^3}}  m^{\alpha} v \, Q_{\beta ,
  \alpha}
(f_{\beta}, f_{\alpha}) dv  =0,
\\ \int_{\mathbb{R^3}}  m^{\alpha} v^2 \, Q_{\alpha ,
  \beta}
(f_{\alpha}, f_{\beta}) dv + \int_{\mathbb{R^3}}  m^{\alpha} v^2 \, Q_{\beta ,
  \alpha}
(f_{\beta}, f_{\alpha}) dv  =0,
\end{eqnarray*}
\end{property}
The function $\mathcal{B}^{\alpha , \beta} (v - v_{\ast},\omega)$ 
is the kernel of the
collision operator $Q_{\alpha,\beta}$. It is a nonegative function
whose form is determined by the molecular interaction between the
species $\alpha$ and $\beta$. Because of the action and reaction
principle, it has the symmetry property $\mathcal{B}^{A , B} =
\mathcal{B}^{B , A}$. More precisely, we consider in this paper the
following type of kernels
\begin{eqnarray*}
\frac{1}{4 \sqrt{2 \pi}} \left(  \frac{d^{\alpha} + d^{\beta}}{2}   \right)^2 |v -
v_{\ast}|^{\beta} b(\theta),
\end{eqnarray*}
 with
\begin{eqnarray*} 0 \leq \beta <
2, \quad b \in L_{+}^{1}([0,2 \pi]), \quad b(\theta) \geq c>0
\quad a.e.
\end{eqnarray*}
for hard forces and 
\begin{eqnarray*} -3 \leq \beta <
0, \quad b \in L_{+}^{1}([0,2 \pi]), \quad b(\theta) \geq c>0
\quad a.e.
\end{eqnarray*}
for soft forces.

As (\cite{ABT}) define the collision frequency as the vector 
$(\nu_A , \nu_B)$, with
\begin{eqnarray*}
 \nu_{\alpha} = \displaystyle \sum_{\beta \in \{ A,B \}}
\int B^{\alpha , \beta} f_{\beta } d\omega dv_{\ast}, \hspace*{1 mm}
\alpha \in \{A,B\}.
\end{eqnarray*}
 On the boundary 
of the domain, the 
two components satisfy different physical properties. Indeed, the $A$ 
component is supposed to be a condensable gas whereas the $B$ component is 
supposed to be a non condensable gas. 

 Hence the boundary condition for the $A$ component is the given indata profile
\begin{eqnarray} \label{bona}
f_{A}(-1,v) = kM_{-}(v) , \xi > 0, \qquad f_{A}(1,v) = kM_{+}(v), \xi < 0,
\end{eqnarray}
for some positive $k$. The boundary condition for the $B$
component is of diffuse reflection type
\begin{eqnarray} \label{bonb}
f_{B}(-1,v) = (\int_{\xi' < 0} | \xi' |  f_{B}(-1,v')
dv')M_{-}(v), \quad \xi>0,
\\ \nonumber f_{B}(1,v) = (\int_{\xi' > 0}  \xi '  f_{B}(1,v') dv')M_{+}(v),\quad \xi<0.
\end{eqnarray}
$M_{+}$ and $M_{-}$ are given normalized Maxwellians
\begin{eqnarray*}
M_{-}(v) = \frac{1}{2\pi T_{-}^{2}} e^{- \frac{|v|^{2}}{2 T_{-}}}
\quad \mbox{and} \quad M_{+}(v) = \frac{1}{2\pi T_{+}^{2}} e^{-
  \frac{|v|^{2}}{2 T_{+}}}.
\end{eqnarray*}
As a theoritical point of view, existence theorem for single component 
gases has been firstly considered. These papers are of interest
because the case of the stationary Boltzmann equation is not covered
by the DiPerna Lions theory established for the time dependant non
linear Boltzmann equation (\cite{DPL}, \cite{CIL}). In (\cite{LA1}), an $L^1$ 
existence theorem is shown for hard and
 soft forces when the distribution function has
a given indatta profile. In
the case of boundary conditions of Maxwell diffuse reflection type, an analogous
theorem is also shown in (\cite{LA2}). In these two papers the solution
are constructed in such a way that they have a given weighted mass.
Let us mention the case of the stationnary Povzner equation in the
case of hard and soft forces which is
investigated in (\cite{P}). The situation of a two component gas has been considered in
(\cite{B1}, \cite{Bw}) when the molecular masses of the two gases are
the same but with different boundary conditions. In these papers, the strategy of the resolution is to use 
that the sum of the distribution of the two components satisfies the
Boltzmann equation for a one component gas. Hence the weak $L^1$ compactness is
firstly obtained for the sum and transmitted to the two distribution
functions. But in the present, case due to the different molecular
masses, the sum of the distribution functions is not the solution of the solution of the Boltzmann
eqution for a single component gas. Therefore the weak $L^1$
compactenss has to be extracted directy on each component. In
(\cite{B2}) the situation of a binary mixture close to a local
equilibrium is investigated. In that 
case the solution of the system is constructed as a Hilbert expansion     
and the rest term is rigorously controled. In \cite{Do} a moment
method is applied in the situation of small Knudsen number to derive a fluid
system. 

As a physical point of view and as a numerical point of view, 
a problem of evaporation condensation for 
a binary mixture far from equilibrium has been considered in
\cite{SAD}. 
The binary mixture
composed of vapor and non condensable gas in contact with an infinite 
plane of condensed vapor. Moreover the non condensable gas is supposed 
to be closed to the condensed vapor. For the numerical simulations 
the authors used a time-dependant BGK
model for a two component gas until a stationary state is reached. The situation of a small Knudsen
number, has also been investigated in (\cite{A}, 
\cite{ATK}, \cite{ATT}, \cite{TAT}) and two types of behaviour 
are pointed out. In a first situation the macroscopic velocity of 
the two gases tends to zero when the Knudsen number tends to zero. 
But the zero order term of the temperature in its Hilbert expansion 
is calculated from the first order term of the macroscopic velocity. 
This means that the macroscopic velocity disappears at the limit 
but keeps an influence on the limit. This is the ghost effect pointed
in \cite{SATSB} for a one component gas and in (\cite{A},\cite{ATK},
\cite{ATT}) for a two component gas. In a second case the B component
becomes negligeable and the macroscopic velocity of the A component 
becomes constent. Moreover the B component accumulates in a thin layer
called Knudsen layer at the boundary where the A component blows. In 
the situation of vapor-vapor mixture ghost-effects have also 
been shown in (\cite{T}). 

 In this paper, weak solutions $(f_{A},f_{B})$ to the stationary problem in the sense of Definition \ref{aa} will be considered.
\begin{definition} \label{aa}
Let $M_{A}$ and $M_{B}$ be given nonnegative real numbers.
$(f_{A},f_{B})$ is a weak solution to the stationary Boltzmann problem with the $\beta$-norms $M_{A}$ and $M_{B}$, if $f_{A}$ and $f_{B}$ $\in$ $L_{loc}^{1}((-1,1)
\times\mathbb{R}^{3})$, $\nu_A , \, \nu_B \,   \in L_{loc}^{1}((-1,1) \times \mathbb{R}^{3}) $, $ \int ( 1 + |v|)^{\beta} f_{A}(x,v) dxdv = M_{A}$, $ \int ( 1 + |v|)^{\beta} f_{B}(x,v)
dxdv = M_{B}$, and there is a constant $k > 0$ such that
for every test function $\varphi \in C_{c}^{1}([-1,1] \times \mathbb{R}^{3}) $  such that $\varphi$ vanishes in a neiborhood of $ \xi = 0$, and on
\\ $ \{ (-1,v) ; \xi < 0 \} \cup  \{ (1,v) ; \xi > 0 \}  $,
\begin{eqnarray} 
 \nonumber \int_{-1}^{1} \int_{\mathbb{R}^{3}} (\xi f_{A} \frac{\partial
\varphi}{\partial x} + Q_{AA}(f_{A},f_{A} )+ Q_{AB}(f_{A},f_{B} ) \varphi)(x,v) dx dv
\\ \nonumber   =  k \int_{\mathbb{R}^{3}, \xi < 0} \xi M_{+}(v) \varphi(1,v) dv
 -  k\int_{\mathbb{R}^{3}, \xi > 0} \xi M_{-}(v) \varphi(-1,v) dv
,
\\ \nonumber \int_{-1}^{1} \int_{\mathbb{R}^{3}} (\xi f_{B}
\frac{\partial \varphi}{\partial x} + Q_{BB}(f_{B}, f_{B})+ 
Q_{BA}(f_{B}, f_{A}) \varphi)(x,v) dx dv ,
\\ \nonumber   =   \int_{\xi ' < 0} | \xi | M_{+}(v) \varphi(1,v) dv (\int_{\xi ' >0} \xi' f_{B}(1,v') dv' )
\\ \nonumber  - \int_{\xi ' > 0}  \xi  M_{-}(v) \varphi(-1,v) dv (\int_{\xi '<0} \xi' f_{B}(-1,v') dv' ).
\end{eqnarray}
\end{definition}
\noindent Renormalized solutions will also been considered. We recall
their definition. Let $g$ be defined for $x>0$ by
\begin{eqnarray*}
g(x) = \ln (1 + x).
\end{eqnarray*}
\begin{definition} \label{an}
Let $M_{A}$ and $M_{B}$ be given nonnegative real numbers.
$(f_{A},f_{B})$ is a renormalized solution to the stationary Boltzmann problem with the $\beta$-norms $M_{A}$ and $M_{B}$, if $f_{A}$ and $f_{B}$ $\in$ $L_{loc}^{1}((-1,1)
\times\mathbb{R}^{3})$, $\nu_A , \, \nu_B \, \in L_{loc}^{1}((-1,1) \times \mathbb{R}^{3}) $, $ \int ( 1 + |v|)^{\beta} f_{A}(x,v) dxdv = M_{A}$, $ \int ( 1 + |v|)^{\beta} f_{B}(x,v)
dxdv = M_{B}$, and there is a constant $k > 0$ such that
for every test function $\varphi \in C_{c}^{1}([-1,1] \times \mathbb{R}^{3}) $  such that $\varphi$ vanishes in a neiborhood of $ \xi = 0$ and on
\\ $ \{ (-1,v) ; \xi < 0 \} \cup  \{ (1,v) ; \xi > 0 \}  $,
\begin{eqnarray} 
 \nonumber \int_{-1}^{1} \int_{\mathbb{R}^{3}} (\xi g(f_{A}) \frac{\partial
\varphi}{\partial x} + \frac{Q_{AA}(f_{A},f_{A} )}{1 + f_{A}} \varphi + 
\frac{Q_{AB}(f_{A},f_{B})}{1 + f_{A}}\, \varphi)(x,v) dx dv
\\ \nonumber   =   \int_{\mathbb{R}^{3}, \xi < 0} \xi g ( k M_{+}(v) ) \varphi(1,v) dv
 -  \int_{\mathbb{R}^{3}, \xi > 0} g (\xi k M_{-}(v)) \varphi(-1,v) dv
,
\\ \nonumber \int_{-1}^{1} \int_{\mathbb{R}^{3}} (\xi g(f_{B}) \frac{\partial \varphi}{\partial x} 
+ \frac{Q_{BB}(f_{B},f_{A} + f_{B})}{1+ f_{B}} \varphi + \frac{Q_{BA}(f_{B},f_{A} )}{1+ f_{B}}\, \varphi)(x,v) dx dv ,
\\ \nonumber   =   \int_{\xi  < 0} \xi  g \big( (\int_{\xi ' >0} \xi' f_{B}(1,v') dv' ) M_{+}(v) \big) \varphi(1,v) dv 
\\ \nonumber  - \int_{\xi  > 0}  \xi g \big(\int_{\xi '<0} \xi' f_{B}(-1,v') dv' ) M_{-}(v) \big) \varphi(-1,v) dv .
\end{eqnarray}
\end{definition}
\noindent The main results of this paper are the following theorems
\begin{theorem} \label{th}
Given $\beta$ with $ 0 \leq \beta <2$, $M_A>0$ and $M_B>0$ there is a
weak solution to the stationary problem with $\beta$-norms equal to
$M_A$ and $M_B$.
\end{theorem}
\begin{theorem} \label{thr}
Given $\beta$ with $  -3 < \beta <0$ $M_A>0$ and $M_B>0$, there is a
renormalized solution to the stationary problem with $\beta$-norms
equal to $M_A$ and $M_B$.
\end{theorem}
The present paper is organized as follows. 
The second and the third section are devoted to the proof of the
theorems \ref{th} and \ref{thr}. In section 2, we perform a fix point step on an approched problem as
in (\cite{LA1}, \cite{LA2}, \cite{B1}, \cite{Bw}). In the last part we
perform a passage to the limit in the sequences of approximation. 
\sectionnew{Approximations with fixed total masses}
Let $r>0, m \in \mathbb{N}^{\ast}, \mu>0, \delta >0, j \in
\mathbb{N}^{\ast}$.
\\ By arguing as in (\cite{LA}), 
we can construct a function, $\chi^{r,m} \in C_{0}^{\infty}$ with range $[0,1]$
invariant under the collision transformations $J_{\alpha , \beta}$,
 for any $\alpha, \beta \in \{ A,B \} $ where
\begin{eqnarray*}
 J_{\alpha , \beta}(v,v_{\ast},\omega) = (v^{(\alpha , \beta ) 
  \prime} 
,v_{\ast}^{(\alpha, \beta ) \prime}, -\omega ), 
\end{eqnarray*}
and under the exchange of $v$ and $v_{\ast}$ and such that
 \begin{eqnarray*}
\chi^{r,m}(v,v_{\ast},\omega) = 1   ,\hspace*{3mm } 
\forall (\alpha , \beta ) \in \{ 
 A,B  \}
\; min(| \xi|, | \xi_{\ast}|,
| \xi^{(\alpha , \beta ) \prime}|,| \xi_{\ast}^{(\alpha ,
  \beta)\prime}| \geq r),
\end{eqnarray*}
and 
\begin{eqnarray*}
 \chi^{r,m}(v,v_{\ast},\omega) = 0   ,\hspace*{3mm } 
\forall (\alpha , \beta ) \in \{ 
 A,B  \}
\;
  max(| \xi|,
| \xi_{\ast}|,| \xi^{\alpha , \beta, \prime}|,| \xi_{\ast}^{\alpha , \beta, \prime}|) \leq r - \frac{1}{m}.
\end{eqnarray*}
The modified collision kernel $\mathcal{B}_{m,n,\mu}^{\alpha, \beta}$ is a positive $C^{\infty}$ function approximating
$min( \mathcal{B}^{\alpha,\beta}, \mu)$, when
 \begin{eqnarray*}
v^{2} + v_{\ast}^{2} < \frac{\sqrt{n}}{2}, \mbox{and } | \frac{v -
v_{\ast}}{|v - v_{\ast} |} \cdot \omega | > \frac{1}{m}, \mbox{
and } | \frac{v - v_{\ast}}{|v - v_{\ast} |} \cdot \omega | <1 -
\frac{1}{m}
\end{eqnarray*}
and such that $ \mathcal{B}_{m,n,\mu}^{\alpha,\beta}(v,v_{\ast},\omega) = 0$,
if
\begin{eqnarray*}
 v^{2} + v_{\ast}^{2} > \sqrt n
\mbox{ or }| \frac{v - v_{\ast}}{|v - v_{\ast} |} \cdot \omega | <
\frac{1}{2m}\mbox{, or }  | \frac{v - v_{\ast}}{|v - v_{\ast} |}
\cdot \omega | > 1 - \frac{1}{2m}.
\end{eqnarray*}
The functions
$\varphi_{l}$ are mollifiers in the $x$-variable defined by
$\varphi_{l}(x) := l \varphi(l x)$, where
\begin{eqnarray*}
\varphi \in C_{0}^{\infty}( \mathbb{R}_{v}^{3}),
\quad support(\varphi) \subset (-1,1), \quad \varphi \geq 0, \quad
\int_{-1}^{1}\varphi(x)dx =1 .
\end{eqnarray*}
For the sake of clarity Theorems \ref{th} and \ref{thr} will be shown
for $M_A = M_B =1$. The passage to general weighted masses is
immediate and we refer to (\cite{LA1}, \cite{LA2}, \cite{B1},
\cite{Bw}). 

 Non negative functions $g_A$, $g_B$ $\in K$ and $\theta \in [0,1]$
being given. By arguing as in \cite{B1}, we can construct $F_{A}$ and
$F_{B}$ 
solutions of
the following boundary value problem
\begin{eqnarray} \label{A}
\nonumber \delta F_{A} +  \xi \frac{\partial}{\partial x}F_{A} =
\int_{ \mathbb{R}_{v_{\ast}}^{3} \times \mathbb{S}^2} \chi^{r,m} \mathcal{B}_{m,n,\mu}^{AA}
 \frac{F_{A}}{1+ \frac{F_{A}}{j}}(x,v') \frac{g_A \ast
  \varphi}{1+\frac{g_A\ast \varphi}{j}}(x,v_{\ast}^{'}) dv_{\ast}
d\omega
\\ \nonumber +\int_{ \mathbb{R}_{v_{\ast}}^{3} \times \mathbb{S}^2} \chi^{r,m} 
\mathcal{B}_{m,n,\mu}^{AB}
 \frac{F_{A}}{1+ \frac{F_{A}}{j}}(x,v') \frac{g_B \ast
  \varphi}{1+\frac{g_B\ast \varphi}{j}}(x,v_{\ast}^{'}) dv_{\ast}
d\omega
\\  \nonumber - F_{A}\int_{ \mathbb{R}_{v_{\ast}}^{3} \times \mathbb{S}^2}
\chi^{r,m} 
\mathcal{B}_{m,n,\mu}^{AA}
\frac{g_A \ast \varphi}{1+\frac{g_A \ast \varphi}{j}}(x,v_{\ast})
dv_{\ast}d\omega
\\  \nonumber  - F_{A}\int_{ \mathbb{R}_{v_{\ast}}^{3} \times \mathbb{S}^2} \chi^{r,m} \mathcal{B}_{m,n,\mu}
\frac{g_B \ast \varphi}{1+\frac{g_B \ast \varphi}{j}}(x,v_{\ast})
dv_{\ast}d\omega, \hspace*{2 mm}
  (x,v) \in (-1,1) \times
\mathbb{R}_{v}^{3},
\\   F_{A}(-1,v)  =  \lambda M_{-}(v), \hspace*{2 mm}  \xi>0, 
\quad F_{A}(1,v)  =  \lambda M_{+}(v)  , \hspace*{2 mm}  \xi<0 , \hspace*{15 mm}
\end{eqnarray}

\noindent and

\begin{eqnarray} \label{B}
\nonumber \delta F_{B} +  \xi \frac{\partial}{\partial x}F_{B} =
\int_{ \mathbb{R}_{v_{\ast}}^{3} \times \mathbb{S}^2} \chi^{r,m} \mathcal{B}_{m,n,\mu}^{BB}
\frac{F_{B}}{1+\frac{F_{B}}{j}}(x,v') \frac{g_B \ast
\varphi}{1+\frac{g_B \ast \varphi}{j}}(x,v_{\ast}^{'})dv_{\ast}
d\omega
\\ \nonumber +
\int_{ \mathbb{R}_{v_{\ast}}^{3} \times \mathbb{S}^2} \chi^{r,m} 
\mathcal{B}_{m,n,\mu}^{AB}
\frac{F_{B}}{1+\frac{F_{B}}{j}}(x,v') \frac{g_A \ast
\varphi}{1+\frac{g_A \ast \varphi}{j}}(x,v_{\ast}^{'})dv_{\ast}
d\omega
\\  \nonumber - F_{B} \int_{ \mathbb{R}_{v_{\ast}}^{3} \times
  \mathbb{S}^2}\chi^{r,m} 
\mathcal{B}_{m,n,\mu}^{BB}
\frac{g_B \ast \varphi}{1+\frac{g_B \ast
\varphi}{j}}(x,v_{\ast})dv_{\ast}d\omega
\\ \nonumber - F_{B} \int_{ \mathbb{R}_{v_{\ast}}^{3} \times 
\mathbb{S}^2}\chi^{r,m} \mathcal{B}_{m,n,\mu}^{BA}
\frac{g_A \ast \varphi}{1+\frac{g_A \ast
\varphi}{j}}(x,v_{\ast})dv_{\ast}d\omega, 
\hspace*{2 mm} (x,v) \in (-1,1)
\times \mathbb{R}_{v}^{3}, 
\\  F_{B}(-1,v) =  \theta \lambda M_{-}(v), \hspace*{2 mm}  \xi>0, 
\hspace*{5 mm}
 F_{B}(1,v) =  (1-\theta)\lambda M_{+}(v)  , \hspace*{2 mm}  \xi<0,
\hspace*{5 mm}
\end{eqnarray}
\noindent  
 as the $L^1$ limit of
sequences. It can also been proven that the equations (\ref{A}) and
(\ref{B}) 
each has a unique solution which is strictly positive. 
\noindent Let
\begin{eqnarray*} f_{A}= \frac{F_{A}}{ \int
min(\mu,(1+|v|)^{\beta})F_{A}(x,v)dx dv } , 
\\ f_{B}=
\frac{F_{B}}{ \int min(\mu,(1+|v|)^{\beta})F_{B}(x,v)dx dv }.
\end{eqnarray*}

\noindent Hence it follows that the functions $f_{A}$ and $f_{B}$ are well defined since
$F_{A}$ and $F_{B}$ strictly positive.

 \noindent Indeed using that 
$\int_{-1}^{1} (\alpha + \nu(x,v))dx \leq 2+2\mu$, it holds that

\begin{eqnarray*} 
F_{A}(x,v) \geq \lambda M_{-}(v) e^{-
\frac{2+2\mu}{\xi}} , \hspace*{2 mm} \xi > 0,
\hspace*{5 mm}  F_{A}(x,v) \geq \lambda M_{+}(v) e^{- \frac{2+2\mu}{| \xi
|}} , \hspace*{2 mm}  \xi <0 .
\end{eqnarray*}

\noindent Analogously, we obtain

\begin{eqnarray*} F_{B}(x,v) \geq \theta \lambda M_{-}(v) e^{-
\frac{2+2 \mu}{\xi}}  , \quad \xi > 0,
\\ F_{B}(x,v) \geq  (1- \theta) \lambda M_{+}(v) e^{- \frac{2+2\mu}{| \xi
|}} ,\quad  \xi <0 .
\end{eqnarray*}

\noindent By taking $ \lambda$ as 

\begin{eqnarray} \label{lam}
\nonumber \lambda = min(\frac{1}{\int_{\xi > 0 }
M_{-}(v)min(\mu,(1+|v|)^{\beta}) e^{- \frac{2+2\mu}{\xi}}dv}; \qquad \qquad \qquad \qquad 
\\ \nonumber \qquad \qquad \frac{1}{\int_{ \xi<0} M_{+}(v) min(\mu,(1+|v|)^{\beta}) e^{-
\frac{2+2\mu}{| \xi |}}dv}),
\end{eqnarray}

\noindent we get

 \begin{eqnarray*} \int min(\mu,(1+|v|)^{\beta})F_{A}(x,v)dx dv
\geq 1
\end{eqnarray*}
\noindent and
\begin{eqnarray*} \int
min(\mu,(1+|v|)^{\beta})F_{B}(x,v)dx dv \geq 1 .
\end{eqnarray*}
\noindent Hence the functions $f_{A}$ and $f_{B}$ are solutions to

\begin{eqnarray} \label{sA}
\nonumber  \delta f_{A} + \xi \frac{ \partial}{\partial x}f_{A} =
\int_{ \mathbb{R}_{v_{\ast}}^{3} \times \mathbb{S}^{2}} \chi^{r,m} 
\mathcal{B}_{m,n,\mu}^{AA}
\frac{f_{A}}{1+\frac{F_A}{j}}(x,v') \frac{g_A \ast \varphi}{1+\frac{g_A
\ast \varphi}{j}}(x,v_{\ast}^{'})dv_{\ast} d\omega
\\ \nonumber  + \int_{ \mathbb{R}_{v_{\ast}}^{3} \times \mathbb{S}^{2}} 
\chi^{r,m} \mathcal{B}_{m,n,\mu}^{AB}
\frac{f_{A}}{1+\frac{F^A}{j}}(x,v') \frac{g_B \ast \varphi}{1+\frac{g_B
\ast \varphi}{j}}(x,v_{\ast}^{'})dv_{\ast} d\omega
\\  \nonumber - f_{A} \int_{ \mathbb{R}_{v_{\ast}}^{3} \times \mathbb{S}^{2}} \chi^{r,m} \mathcal{B}_{m,n,\mu}^{AA}
\frac{g_A \ast \varphi}{1+\frac{g_A \ast
    \varphi}{j}}(x,v_{\ast})dv_{\ast} d\omega
\\ \nonumber - f_{A} \int_{ \mathbb{R}_{v_{\ast}}^{3} \times \mathbb{S}^{2}} 
\chi^{r,m} \mathcal{B}_{m,n,\mu}^{AB}
\frac{g_B \ast \varphi}{1+\frac{g_B \ast
    \varphi}{j}}(x,v_{\ast})dv_{\ast} d\omega,
\hspace*{3 mm} (x,v) \in
(-1,1) \times \mathbb{R}_{v}^{3}, 
\\  \nonumber f_{A}(-1,v) = \frac{\lambda}{ \int
min(\mu,(1+|v|)^{\beta})F_{A}(x,v)dx dv} M_{-}(v)  , \hspace*{2 mm} \xi>0,
\\  \nonumber f_{A}(1,v) =  \frac{\lambda}{ \int
min(\mu,(1+|v|)^{\beta})F_{A}(x,v)dx dv} M_{+}(v) , \hspace*{2 mm} \xi <
0,
\\ 
\end{eqnarray}

\noindent and

\begin{eqnarray} \label{sB}
\nonumber \delta f_{B} +\xi \frac{\partial}{\partial x}f_{B} =
 \int_{ \mathbb{R}_{v_{\ast}}^{3} \times \mathbb{S}^{2}} \chi^{r,m}
 \mathcal{B}_{m,n,\mu}^{BB}(v,v_{\ast},\omega) \frac{f_{B}}{1+\frac{F_B}{j}}(x,v')
 \frac{g_B \ast \varphi}{1+\frac{g_B \ast
     \varphi}{j}}(x,v_{\ast}^{'})dv_{\ast} d\omega
\\ \nonumber + 
 \int_{ \mathbb{R}_{v_{\ast}}^{3} \times \mathbb{S}^{2}} \chi^{r,m}
 \mathcal{B}_{m,n,\mu}^{BA}(v,v_{\ast},\omega) \frac{f_{B}}{1+\frac{F_B}{j}}(x,v')
 \frac{g_B\ast \varphi}{1+\frac{g_B \ast
     \varphi}{j}}(x,v_{\ast}^{'})dv_{\ast} d\omega
\\  \nonumber - f_{B}(x,v) \int_{ \mathbb{R}_{v_{\ast}}^{3} \times
  \mathbb{S}^{2}}\chi^{r,m}
 \mathcal{B}_{m,n,\mu}^{BB}
\frac{g_B \ast \varphi}{1+\frac{g_B \ast
\varphi}{j}}(x,v_{\ast})dv_{\ast} d\omega
\\  \nonumber- f_{B}(x,v) \int_{ \mathbb{R}_{v_{\ast}}^{3} 
\times \mathbb{S}^{2}}\chi^{r,m} \mathcal{B}_{m,n,\mu}^{BA}
\frac{g_A \ast \varphi}{1+\frac{g_A \ast
\varphi}{j}}(x,v_{\ast})dv_{\ast} d\omega,  \hspace*{2 mm} (x,v) \in (-1,1)
\times \mathbb{R}_{v}^{3}, \quad
\\ \nonumber f_{B}(-1,v)= \frac{\lambda}{ \int min(\mu,(1+|v|)^{\beta})F_{B}(x,v)dx dv} \theta M_{-}(v) , \quad \xi>0,
\\ \nonumber f_{B}(1,v) =  \frac{\lambda}{ \int
min(\mu,(1+|v|)^{\beta})F_{B}(x,v)dx dv} (1- \theta) M_{+}(v) ,
\quad \xi < 0.
\\
\end{eqnarray}

\noindent In order to use a fixed-point theorem, consider the closed and convex subset of $L_{+}^{1}([-1,1] \times \mathbb{R}_{v}^{3}) $,
\begin{eqnarray*} 
K =\{ f \in L_{+}^{1}([-1,1] \times \mathbb{R}_{v}^{3}), \quad
\int_{[-1,1] \times \mathbb{R}_{v}^{3}} min( \mu, (1+|v |)^{\beta}) f(x,v) dxdv=1 \}.
\end{eqnarray*} 
The fixed-point argument will now be used in order to
solve (\ref{sA}, \ref{sB}) with $g_A = f_{A}$ and $g_B=f_{B}$.

 Define $T$ on $K \times K \times [0,1]$ by $ T(g_A,g_B, \theta) = (f_{A} , f_{B} , \tilde{\theta}) $ with
\begin{eqnarray} \label{t}
\tilde{\theta} = \frac{ \int_{ \xi <0}| \xi | f_{B}(-1,v)dv}{\int_{
  \xi <0} |\xi | f_{B}(-1,v)dv + \int_{ \xi >0} \xi f_{B}(1,v)dv}
\end{eqnarray}
where $(f_{A},f_{B})$ is solution to (\ref{sA}, \ref{sB}).

By reasonning as in \cite{B1}, it can be shown that the 
map $T$ is continous from $K \times K \times [0,1]$ into itself. So
from the Schauder fixed point theorem there is $(f_A,f_B,\theta)$ 
such that 
\begin{eqnarray*}
 f_A = g_A,\hspace*{4 mm}
 f_B = g_B ,
\hspace*{4 mm}  \theta = \frac{ \int_{ \xi<0}| \xi |f_{B}(-1,v)dv }{ \int_{\xi>0}\xi
   f_{B}(1,v)dv+{ \int_{ \xi<0}| \xi |f_{B}(-1,v)dv}} 
\end{eqnarray*}
 that satisfy
\begin{eqnarray} \label{qa}
\nonumber \delta f_{A} +\xi \frac{\partial}{\partial x}f_{A} =
\int_{ \mathbb{R}_{v_{\ast}}^{3} \times \mathbb{S}^{2}} \chi^{r,m} 
 \mathcal{B}_{m,n,\mu}^{AA}
\frac{f_{A}}{1+\frac{F_A}{j}}(x,v') \frac{f_A \ast
\varphi_{l}}{1+\frac{f_A \ast
\varphi_{l}}{j}}(x,v_{\ast}^{'})dv_{\ast} d\omega \qquad
\\ \nonumber + \int_{ \mathbb{R}_{v_{\ast}}^{3} \times \mathbb{S}^{2}} \chi^{r,m} 
\mathcal{B}_{m,n,\mu}^{AB}
\frac{f_{A}}{1+\frac{F_A}{j}}(x,v') \frac{f_B \ast
\varphi_{l}}{1+\frac{f_B \ast
\varphi_{l}}{j}}(x,v_{\ast}^{'})dv_{\ast} d\omega \qquad
\\   \nonumber - f_{A}  \int_{ \mathbb{R}_{v_{\ast}}^{3} \times
  \mathbb{S}^{2}}\chi^{r,m} 
\mathcal{B}_{m,n,\mu}^{AA}
\frac{f_A \ast \varphi_{l}}{1+\frac{f_A \ast
\varphi_{l}}{j}}(x,v_{\ast})dv_{\ast} d\omega
 \\ \nonumber - f_{A}  \int_{ \mathbb{R}_{v_{\ast}}^{3} \times \mathbb{S}^{2}}\chi^{r,m} \mathcal{B}_{m,n,\mu}^{AB}
\frac{f_B \ast \varphi_{l}}{1+\frac{f_B \ast
\varphi_{l}}{j}}(x,v_{\ast})dv_{\ast} d\omega   ,
\hspace*{2 mm}
 (x,v) \in
(-1,1) \times \mathbb{R}_{v}^{3}, \quad
\\ \nonumber f_{A}(-1,v) = k_{A} M_{-}(v), \hspace*{2 mm} \xi>0, \quad
f_{A}(1,v) =   k_{A} M_{+}(v), \hspace*{2 mm} \xi <
0 \hspace*{4 mm}
\\
\end{eqnarray}
with
\begin{eqnarray*}
 k_{A} = \frac{\lambda}{\int min(\mu , (1 + |v|)^{\beta})
  F_{A}(x,v)dx dv}
\end{eqnarray*}
\noindent and
\begin{eqnarray} \label{qb}
\nonumber \delta f_{B} +\xi \frac{\partial}{\partial x}f_{B} =
\int_{ \mathbb{R}_{v_{\ast}}^{3} \times \mathbb{S}^{2}} \chi^{r,m} 
\mathcal{B}_{m,n,\mu}^{BB}
\frac{f_{B}}{1+\frac{F^B}{j}}(x,v') \frac{f_B \ast
\varphi_{l}}{1+\frac{f_B \ast
\varphi_{l}}{j}}(x,v_{\ast}^{'})dv_{\ast} d\omega \qquad
\\ \nonumber + \int_{ \mathbb{R}_{v_{\ast}}^{3} \times \mathbb{S}^{2}}
\chi^{r,m} 
\mathcal{B}_{m,n,\mu}^{BA}
\frac{f_{B}}{1+\frac{F_B}{j}}(x,v') \frac{f_A \ast
\varphi_{l}}{1+\frac{f_A \ast
\varphi_{l}}{j}}(x,v_{\ast}^{'})dv_{\ast} d\omega \qquad
\\ \nonumber   - f_{B} \int_{ \mathbb{R}_{v_{\ast}}^{3} 
\times \mathbb{S}^{2}}\chi^{r,m} \mathcal{B}_{m,n,\mu}^{BB}
\frac{f_B \ast \varphi_{l}}{1+\frac{f_B \ast
    \varphi_{l}}{j}}(x,v_{\ast})dv_{\ast} d\omega 
\\ \nonumber - f_{B} \int_{ \mathbb{R}_{v_{\ast}}^{3} \times
  \mathbb{S}^{2}}\chi^{r,m} 
 \mathcal{B}_{m,n,\mu}^{BA}
\frac{f_A \ast \varphi_{l}}{1+\frac{f_A \ast
    \varphi_{l}}{j}}(x,v_{\ast})dv_{\ast} d\omega , \hspace*{3 mm}
   (x,v) \in
(-1,1) \times \mathbb{R}_{v}^{3}, 
\\ \nonumber f_{B}(-1,v)  =  \lambda '  ( \frac{ \int_{ \xi<0}| \xi |f_{B}(-1,v)dv }{ \int_{\xi>0}\xi f_{B}(1,v)dv+{ \int_{ \xi<0}| \xi
|f_{B}(-1,v)dv}})M_{-}(v), \hspace*{2 mm} \xi > 0,
\\ \nonumber f_{B}(1,v)  = \lambda '  (\frac{ \int_{ \xi>0}| \xi
  |f_{B}(1,v)dv }{ \int_{\xi>0}\xi f_{B}(1,v)dv+{ \int_{ \xi<0}| \xi
    |f_{B}(-1,v)dv}} ) M_{+}(v) , \hspace*{2 mm}
\xi < 0 , \quad
\\ 
\end{eqnarray}
with
\begin{eqnarray*}
 \lambda ' = \frac{\lambda}{ \int
min(\mu,(1+|v|)^{\beta})F_{B}(x,v)dx dv}.
\end{eqnarray*}
\section{The slab solution for $-3 < \beta \leq 0 $ and $0 \leq \beta
  <2$.}
This section is devoted to the passage to the limit in 
(\ref{qa}, \ref{qb}). It is performed in two times. In the first one 
the solutions of the approached problem are written in their
exponential form and averaging lemmas are applied.
The second passage to the limit corresponds to the passage to the limit in (\ref{la},
\ref{lb}). One crucial point is to get an entropy estimate on 
$(f_A^j,f_B^j)$ in order to extract compactness. In
(\cite{B1}), this control is obtained from a bound on the entropy of 
$f^j = f_A^j + f_B^j$ by using that $f^j$ satisfy the Boltzmann
equation for a single component gas. But in the present paper, 
due to the difference of the molecular masses, this property is 
not satisfied.   

Keeping, $l$, $j$, $r$, $m$, $\mu$ fixed, denote
$f^{j,\delta,l,r,m,\mu}$ by $f^{\delta}$ each distribution function and study the passage to the
limit when 
$\delta$ tends to 0.
Writing the equations (\ref{qa}, \ref{qb}) in the exponential form and using the averaging lemmas together with a convolution with a mollifier
(\cite{LA2},\cite{L}) give that $f_{A}^{\delta}$ and $F_A^{\delta}$ are
strongly compact in $L^{1}([-1,1] \times \mathbb{R}_{v}^{3})$. Denote by
 $f_{A}$ and $F_A$ the respective limits of $f_{A}^{\delta}$ and $F_A^{\delta}$. The
 passage to the limit when $\delta$ tends to 0 in the equation
 (\ref{qa}) yields
\begin{eqnarray} \label{la}
\nonumber  \xi \frac{\partial}{\partial x}f_{A}  =  \int_{
\mathbb{R}_{v_{\ast}}^{3} \times \mathbb{S}^2} \chi^{r,m}
\mathcal{B}_{m,n,\mu}^{AA}
\frac{f_{A}}{1+\frac{F^A}{j}}(x,v') \frac{f_A \ast
\varphi_{l} }{1+\frac{f_A \ast \varphi_{l}
}{j}}(x,v_{\ast}^{'})dv_{\ast} d\omega \qquad
\\ \nonumber \int_{
\mathbb{R}_{v_{\ast}}^{3} \times \mathbb{S}^2} \chi^{r,m}
\mathcal{B}_{m,n,\mu}^{AB} \frac{f_{A}}{1+\frac{F^A}{j}}(x,v') \frac{f_B \ast
\varphi_{l} }{1+\frac{f_B \ast \varphi_{l}
}{j}}(x,v_{\ast}^{'})dv_{\ast} d\omega \qquad
\\ \nonumber    -  f_{A} \int_{ \mathbb{R}_{v_{\ast}}^{3} \times \mathbb{S}^2}\chi^{r,m} \mathcal{B}_{m,n,\mu}^{AA}
\frac{f_A \ast \varphi_{l}}{1+\frac{f_A \ast
\varphi_{l}}{j}}(x,v_{\ast})dv_{\ast} d\omega, 
\\ \nonumber -  f_{A} \int_{ \mathbb{R}_{v_{\ast}}^{3} \times
  \mathbb{S}^2}
\chi^{r,m} \mathcal{B}_{m,n,\mu}^{AB}
\frac{f_B \ast \varphi_{l}}{1+\frac{f_B \ast
\varphi_{l}}{j}}(x,v_{\ast})dv_{\ast} d\omega, 
\hspace*{2 mm}
 (x,v) \in
(-1,1) \times \mathbb{R}_{v}^{3}, \quad
\\ \nonumber f_{A}(-1,v)  =  \frac{\lambda}{ \int
  min(\mu,(1+|v|)^{\beta})F_{A}(x,v)dx dv} M_{-}(v) , \hspace*{2 mm}
\xi>0 ,\qquad
\\ \nonumber f_{A}(1,v)  =  \frac{\lambda}{ \int
min(\mu,(1+|v|)^{\beta})F_{A}(x,v)dx dv} M_{+}(v) , \quad
\xi < 0, 
\\
\end{eqnarray}

with
\begin{eqnarray*}
  \int min(\mu,(1+| v |)^{\beta})f_{A}^j(x,v)dxdv=1.
\end{eqnarray*}

\noindent For the same reasons, the limit $f_{B}$ of
$f_{B}^{\delta}$ satisfies
\begin{eqnarray} \label{lb}
\nonumber \xi \frac{\partial}{\partial x}f_{B} = \int_{
\mathbb{R}_{v_{\ast}}^{3} \times \mathbb{S}^{2}} \chi^{r,m} 
\mathcal{B}_{m,n,\mu}^{BB}
\frac{f_{B}}{1+\frac{F^B}{j}}(x,v') \frac{f_B \ast
\varphi_{l}}{1+\frac{f_B \ast
\varphi_{l}}{j}}(x,v_{\ast}^{'})dv_{\ast} d\omega 
\\ \nonumber + \int_{
\mathbb{R}_{v_{\ast}}^{3} \times \mathbb{S}^{2}} \chi^{r,m} \mathcal{B}_{m,n,\mu}^{BB}
\frac{f_{B}}{1+\frac{F^B}{j}}(x,v') \frac{f_A \ast
\varphi_{l}}{1+\frac{f_A \ast
\varphi_{l}}{j}}(x,v_{\ast}^{'})dv_{\ast} d\omega 
\\   \nonumber -  f_{B} \int_{ \mathbb{R}_{v_{\ast}}^{3} \times
  \mathbb{S}^{2}}\chi^{r,m} 
\mathcal{B}_{m,n,\mu}^{BB}
\frac{f_B \ast \varphi_{l}}{1+\frac{f_B \ast
\varphi_{l}}{j}}(x,v_{\ast})dv_{\ast}d\omega
\\ \nonumber -  f_{B} \int_{ \mathbb{R}_{v_{\ast}}^{3} \times
  \mathbb{S}^{2}}\chi^{r,m} 
\mathcal{B}_{m,n,\mu}^{BA}
\frac{f_A \ast \varphi_{l}}{1+\frac{f_A \ast
\varphi_{l}}{j}}(x,v_{\ast})dv_{\ast}d\omega, \hspace*{2 mm}
 (x,v) \in
(-1,1) \times \mathbb{R}_{v}^{3}, \qquad
\\ \nonumber f_{B}(-1,v) = \sigma(-1) \lambda ' M_{-}(v) , \hspace*{2 mm} \xi>0,
\quad  f_{B}(1,v) = \sigma(1) \lambda ' M_{+}(v)  , \hspace*{2 mm}  \xi<0,
\\
\end{eqnarray}

\noindent with
\begin{eqnarray*}
 \int min( \mu , (1 + | v | )^{ \beta}) f_{B}(x,v)dxdv=1,
\end{eqnarray*}

\noindent where
\begin{eqnarray}
 \nonumber \sigma(-1) &=& \frac{ \int_{ \xi<0}| \xi
|f_{B}(-1,v)dv }{
  \int_{\xi>0}\xi f_{B}(1,v)dv+{ \int_{ \xi<0}| \xi |f_{B}(-1,v)dv}},
\\ \nonumber \sigma^j(1) &=& \frac{ \int_{ \xi>0} \xi f_{B}(1,v)dv }{
    \int_{\xi>0}\xi f_{B}(1,v)dv+{ \int_{ \xi<0}| \xi |f_{B}(-1,v)dv}}
\end{eqnarray}
\noindent and
 \begin{eqnarray*}
\lambda ' = \frac{\lambda}{\int
  min(\mu,(1+|v|)^{\beta})F_{B}^j(x,v)dxdv} .
\end{eqnarray*}
Mutltiply (\ref{la}) by $ \log(\frac{f_A^j}{1 + \frac{f_A^j}{j}} ) $ and 
(\ref{lb}) by $ \log(\frac{f_B^j}{1 + \frac{f_B^j}{j}} ) $ and the two
equations leads to according to (\cite{LA1}, \cite{ABT}, \cite{Do}),
\begin{eqnarray*}
\int_{\mathbb{R}^3} \xi \, \left(   
f_A^j \log (f_A^j)(1,v) - j (1 + \frac{f_A^j}{j} )
\log (1 + \frac{f_A^j}{j})(1,v)
\right)
\\ -\int_{\mathbb{R}^3} \xi \, \left(   
f_A^j \log (f_A^j)(-1,v) - j (1+ \frac{f_A^j}{j} )
\log (1 + \frac{f_A^j}{j})(-1,v)
\right)
\\ + \int_{\mathbb{R}^3} \xi \, \left(   
f_B^j \log (f_B^j)(1,v) - j (1 + \frac{f_B^j}{j}) 
\log (1 + \frac{f_B^j}{j})(1,v)
\right)
\\ -\int_{\mathbb{R}^3} \xi \, \left(   
f_B^j \log (f_B^j)(1,v) - j (1 + \frac{f_B^j}{j}) 
\log (1 + \frac{f_B^j}{j})(1,v)
\right)
\\ = -\frac{1}{4}I_{AA}^j(f_A^j,f_A^j) 
-\frac{1}{2} I_{AB}^j(f_A^j,f_B^j) -\frac{1}{4} I_{BB}^j(f_B^j,f_B^j) 
\\ +  \int
\chi^{r,m} 
\mathcal{B}_{m,n,\mu}^{AA}
\frac{f_A^{j\prime}( f_A^{j\prime} - F_A^{j\prime})}{j(1+  F_A^{j\prime})
(1+  f_A^{j\prime})} \frac{f_{A\ast}^\prime}{1 +
\frac{f_{A\ast}^{j\prime}}{j} }
\log \frac{f_{A}^j}{1 +
\frac{f_{A}^j}{j} }
\\ +  \int
\chi^{r,m} 
\mathcal{B}_{m,n,\mu}^{AB}
\frac{f_A^{j\prime}( f_A^{j\prime} - F_A^{j\prime})}{j(1+  F_A^{j\prime})
(1+  f_A^{j\prime})} \frac{f_{B\ast}^{j\prime}}{1 +
\frac{f_{B\ast}^\prime}{j} }
\log \frac{f_{A}^j}{1 +
\frac{f_{A}^j}{j} }
\\ - \int
\chi^{r,m} 
 \frac{f_A^{j2}}{ j(1+\frac{f_A^j}{j}  ) } 
 \log \frac{f_{A}^j}{1 +
\frac{f_{A}^j}{j} } \left(           
\mathcal{B}_{m,n,\mu}^{AA}
\frac{f_{A\ast}^j}{ (1+\frac{f_{A\ast}^j}{j}  ) } 
  \mathcal{B}_{m,n,\mu}^{AB}   \frac{f_{B\ast}^j}{ (1+\frac{f_{B\ast}^j}{j}  ) }     \right)                        
\\ +  \int
\chi^{r,m} 
\mathcal{B}_{m,n,\mu}^{BB}
\frac{f_B^{j\prime}( f_B^{j\prime} - F_B^{j\prime})}{j(1+  F_B^{j\prime})
(1+  f_B^{j\prime})} \frac{f_{B\ast}^{j\prime}}{1 +
\frac{f_{B\ast}^\prime}{j} }
\log \frac{f_{B}^j}{1 +
\frac{f_{B}^j}{j} }
\\ +  \int
\chi^{r,m} 
\mathcal{B}_{m,n,\mu}^{BA}
\frac{f_B^{\prime}( f_B^{\prime} - F_B^{\prime})}{j(1+  F_B^{\prime})
(1+  f_B^{\prime})} \frac{f_{A\ast}^\prime}{1 +
\frac{f_{A\ast}^\prime}{j} }
\log \frac{f_{B}}{1 +
\frac{f_{B}}{j} }
\\ - \int
\chi^{r,m} 
\mathcal{B}_{m,n,\mu}^{BB} \frac{f_B^{j2}}{ j(1+\frac{f_B^j}{j}  ) } 
\frac{f_{B\ast}^j}{ (1+\frac{f_{B\ast}^j}{j}  ) } 
 \log \frac{f_{B}^j}{1 +
\frac{f_{B}^j}{j} }
\\ - \int
\chi^{r,m} 
\mathcal{B}_{m,n,\mu}^{BA} \frac{f_B^{j2}}{ j(1+\frac{f_B^j}{j}  ) } 
\frac{f_{A\ast}^j}{ (1+\frac{f_{A\ast}^j}{j}  ) } 
 \log \frac{f_{B}^j}{1 +
\frac{f_{B}^j}{j} }
\end{eqnarray*}
with 
\begin{eqnarray*}
I_{AA}^j (f_A^j,f_A^j)=  \int
\chi^{r,m} 
\mathcal{B}_{m,n,\mu}^{AA} 
\left(  
\frac{f_{A}^{j\prime}}{1+\frac{f_A^{j\prime}}{j}}
\frac{f_{A \ast}^{j\prime}}{1+\frac{f_{A \ast}^{j\prime}}{j}} - 
 \frac{f_{A}^j}{1+\frac{f_{A}^j}{j}}
\frac{f_{A\ast}^j}{1+\frac{f_{A\ast}^j}{j}}
\right)  
 \\ \log 
\left(  \frac{\frac{f_A^{\prime}}{1+\frac{f_A^{\prime}}{j}}
    \frac{f_{A \ast}^{j\prime}}{1+\frac{f_{A\ast}^{j\prime}}{j}} }{ \frac{f_A^j}{1+\frac{f_A^j}{j}}
\frac{f_{A\ast}^j}{1+\frac{f_{A\ast}^j}{j}}       }             \right)dx
dv dv_{\ast} d\omega ,
\end{eqnarray*}
\begin{eqnarray*}
I_{BB}^j (f_B^j,f_B^j)=  \int
\chi^{r,m} 
\mathcal{B}_{m,n,\mu}^{BB} 
\left(  
\frac{f_{B}^{\prime}}{1+\frac{f_B^{\prime}}{j}}
\frac{f_{B \ast}^{\prime}}{1+\frac{f_{B \ast}^{\prime}}{j}} - 
 \frac{f_{B}^j}{1+\frac{f_{B}^j}{j}}
\frac{f_{B\ast}^j}{1+\frac{f_{B\ast}^j}{j}}
\right)  
 \\ \log 
\left(  \frac{\frac{f_B^{j\prime}}{1+\frac{f_B^{j\prime}}{j}}
    \frac{f_{B\ast}^{j\prime}}{1+\frac{f_{B\ast}^{j\prime}}{j}} }{ \frac{f_B^j}{1+\frac{f_B^j}{j}}
\frac{f_{B\ast}^j}{1+\frac{f_{B\ast}^j}{j}}       }             \right)dx dv dv_{\ast} d\omega,
\end{eqnarray*}
\begin{eqnarray*}
I_{AB}^j (f_A^j,f_B^j)=  \int
\chi^{r,m} 
\mathcal{B}_{m,n,\mu}^{AB} 
\left(  
\frac{f_{A}^{j\prime}}{1+\frac{f_A^{j\prime}}{j}}
\frac{f_{B \ast}^{j\prime}}{1+\frac{f_{B \ast}^{j\prime}}{j}} - 
 \frac{f_{A}^j}{1+\frac{f_{A}^j}{j}}
\frac{f_{B\ast}^j}{1+\frac{f_{B\ast}^j}{j}}
\right)  
 \\ \log 
\left(  \frac{\frac{f_A^{j\prime}}{1+\frac{f_A^{j\prime}}{j}}
    \frac{f_{B\ast}^{j\prime}}{1+\frac{f_{B\ast}^{j\prime}}{j}} }{ \frac{f_A}{1+\frac{f_A^j}{j}}
\frac{f_{B\ast}^j}{1+\frac{f_{B\ast}^j}{j}}       }             \right)dx dv dv_{\ast} d\omega.
\end{eqnarray*}
According to \cite{ABT}, we have 
$I_{AA}^j(f_A^j,f_A^j) \geq 0$, $I_{AB}^j(f_A^j,f_B^j) \geq 0$ 
$I_{BB}^j(f_B^j,f_B^j) \geq 0$ and by reasonning as in (\cite{LA1}), we can prove that the terms 
\begin{eqnarray*}
 \int
\chi^{r,m} 
\mathcal{B}_{m,n,\mu}^{\alpha \beta} \frac{f_\alpha^2}{ j(1+\frac{f_\alpha}{j}  ) } 
\frac{f_{\beta\ast}}{ j(1+\frac{f_{\beta \ast}}{j}  ) } 
 \log \frac{f_{ \alpha}}{1 +
\frac{f_{\alpha}}{j} },
\\ \int
\chi^{r,m} 
\mathcal{B}_{m,n,\mu}^{\alpha , \beta}
\frac{f_\alpha^{\prime}( f_\alpha^{\prime} - F_\alpha^{\prime})}{j(1+  F_\alpha^{\prime})
(1+  f_\alpha^{\prime})} \frac{f_{\beta\ast}^\prime}{1 +
\frac{f_{\beta\ast}^\prime}{j} }
\log \frac{f_{\alpha}}{1 +
\frac{f_{\alpha}}{j} }
\end{eqnarray*}
are controled uniformly in $j$.
Therefore
\begin{eqnarray*}
 \int_{\mathbb{R}^3} \xi \, \left(   
f_A^j \log (f_A^j)(1,v) - j (1 + \frac{f_A^j}{j} )
\log (1 + \frac{f_A^j}{j})(1,v)
\right)
\\ -\int_{\mathbb{R}^3} \xi \, \left(   
f_A^j \log (f_A^j)(-1,v) - j (1+ \frac{f_A^j}{j} )
\log (1 + \frac{f_A^j}{j})(-1,v)
\right)
\\ + \int_{\mathbb{R}^3} \xi \, \left(   
f_B^j \log (f_B^j)(1,v) - j (1 + \frac{f_B^j}{j}) 
\log (1 + \frac{f_B^j}{j})(1,v)
\right)
\\ -\int_{\mathbb{R}^3} \xi \, \left(   
f_B^j \log (f_B^j)(1,v) - j (1 + \frac{f_B^j}{j}) 
\log (1 + \frac{f_B^j}{j})(-1,v)
\right) \leq c
\end{eqnarray*}
So as in (\cite{LA1}, \cite{LA2}), it follows that $f_A^j$ and $f_B^j$ are weakly compact in $L^1$. 
\begin{remark}
Contrarily to (\cite{B1}, \cite{Bw}), the weak compactness 
of $f_A^j$ and $f_B^j$ is directly obtained. In (\cite{B1}, \cite{Bw}), the author shows 
that the sum $f^j = f_A^j + f_B^j$ is weakly compact in $L^1$ by using 
that $f^j$ satisfies the Boltzmann equation for a single component 
gas. In the present paper, the 2 components having different molecular 
masses, $f^j$ is not solution of the Boltzmann equation for a one 
component gas.     
\end{remark}
Let $Q_{\alpha , \beta}^{j-}$ and  $Q_{\alpha , \beta}^{j+}$ be
defined by 
\begin{eqnarray*}
Q_{\alpha, \beta}^{j-}(f_{\alpha}^j, f_{\beta}^j) = 
f_{\alpha}^j(x,v) \int_{\mathbb{R}^3 \times \mathbb{S}^2} \chi^{r,m}
\mathcal{B}_{m,n,\mu } \frac{f_{\beta}^j}{1 +\frac{f_{\beta}^j}{j}
}(x,v_{\ast})
dv_{\ast} d\omega,
\\ Q_{\alpha, \beta}^{j+}(f_{\alpha}^j, f_{\beta}^j) = 
\int_{\mathbb{R}^3 \times \mathbb{S}^2} \chi^{r,m}
\mathcal{B}_{m,n,\mu } \frac{f_{\alpha}^j}{1 +\frac{f_{\alpha}^j}{j}
}(x,v^{\prime}) \frac{f_{\beta}^j}{1 +\frac{f_{\beta}^j}{j}
}(x,v_{\ast}^{\prime})
dv_{\ast} d\omega.
\end{eqnarray*}
\begin{remark} \label{entr}
The quantity $\frac{1}{4}I_{AA}^j(f_A^j,f_A^j) 
+\frac{1}{2} I_{AB}^j(f_A^j,f_B^j) +\frac{1}{4} I_{BB}^j(f_B^j,f_B^j)
$ is a generalization of the entropy production term used in (\cite{LA1}).
\end{remark}
In order to pass to the limit in (\ref{la}, \ref{lb}) weak compactness
is required on the terms $Q_{\alpha, \beta}^{j-}$ and $Q_{\alpha ,
  \beta}^{j +}$. The inequalities 
\begin{eqnarray*}
Q_{\alpha,\beta}^{j-}(f_\alpha^j,f_\beta^j) 
\leq \int |v-v_{\ast}|^{\beta} f_{\beta}^j
dv_{\ast} d\omega , \; \{\alpha , \beta \} \in \{A,B \},
\end{eqnarray*} 
gives that $Q_{\alpha,\beta}^{j-}$ is weakly compact in $L^1$. 
 By arguing as in \cite{Bw}, we can show that 
\begin{eqnarray} \label{eq1}
\nonumber Q_{A,A}^{j+}(f_A^j,f_A^j)+
Q_{A,B}^{j+}(f_A^j,f_B^j) +
Q_{B,A}^{j+}(f_B^j,f_A^j) +
Q_{B,A}^{j+}(f_B^j,f_A^j)
\\ \nonumber \leq K \left( Q_{A,A}^{j-}(f_A^j,f_A^j) +
Q_{A,B}^{j-}(f_A^j,f_B^j)+
Q_{B,A}^{j-}(f_B^j,f_A^j)+
Q_{B,A}^{j-}(f_B^j,f_A^j) \right)
\\ + \frac{1}{\ln K} \left(
I_{AA}( f_A^j,f_A^j ) + I_{BB}( f_B^j,f_B^j )+ I_{BA}( f_B^j,f_A^j )
\right). \hspace*{5 mm}
\end{eqnarray}
From the weak compactness of $Q_{\alpha,\beta}^{j-}$ for $\{ \alpha ,
\beta \} \in \{A,B \}$ 
and the boundeness from above of
\begin{eqnarray*}
I_{AA}( f_A^j,f_A^j ) + I_{BB}( f_B^j,f_B^j )+ I_{BA}( f_B^j,f_A^j ),
\end{eqnarray*} 
the gain terms 
$Q_{\alpha,\beta}^{j+}$ are weakly compact in $L^1$. Hence by arguing as in (\cite{LA1}, 
\cite{LA2}) we can pass to the limit in the equations 
(\ref{la}, \ref{lb}) and we get that there is $(f_A^{r,\mu},f_B^{r,\mu})$ solution to 
\begin{eqnarray} \label{ra}
\nonumber \xi \frac{\partial}{\partial x}f_{A}^{r,\mu} = \int_{
\mathbb{R}_{v}^{3} \times \mathbb{S}^{2}} \chi^{r}
\mathcal{B}_{\mu}^{AA}(v-v_{\ast},\omega)f_{A}^{r,\mu}(x,v')
f_A^{r,\mu}(x,v_{\ast}^{'})dv_{\ast} d\omega 
\\ \nonumber + \int_{
\mathbb{R}_{v}^{3} \times \mathbb{S}^{2}} \chi^{r}
\mathcal{B}_{\mu}^{AB}(v-v_{\ast},\omega)f_{A}^{r,\mu}(x,v')
f_B^{r,\mu}(x,v_{\ast}^{'})dv_{\ast} d\omega 
\\  \nonumber   - f_{A}^{r,\mu} \int_{ \mathbb{R}_{v_{\ast}}^{3} \times \mathbb{S}^{2}}\chi^{r} \mathcal{B}_{\mu}(v-v_{\ast},\omega)
f_A^{r,\mu}(x,v_{\ast})dv_{\ast}d\omega,
  \\  \nonumber   - f_{A}^{r,\mu} \int_{ \mathbb{R}_{v_{\ast}}^{3} \times \mathbb{S}^{2}}\chi^{r} \mathcal{B}_{\mu}^{AB}(v-v_{\ast},\omega)
f_B^{r,\mu}(x,v_{\ast})dv_{\ast}d\omega ,\quad (x,v) \in (-1,1)
\times \mathbb{R}_{v}^{3}, \quad
\\ f_{A}^{r,\mu}(-1,v) = k_{A}  M_{-}(v), \hspace*{2 mm} \xi>0,
\hspace*{5 mm}
f_{A}^{r,\mu}(1,v)=  k_{A} M_{+}(v)  ,  \hspace*{2 mm}  \xi<0, \hspace*{7 mm}
\end{eqnarray}
\noindent with
 \begin{eqnarray*}
\int min( \mu , (1 + | v | )^{ \beta}) f_{A}^{r,\mu}(x,v)dx
dv=1,
\end{eqnarray*}
 \noindent where $k_{A}$ is defined in the equation (\ref{qa}) before
 passing to the limit.
\begin{eqnarray} \label{rb}
\nonumber \xi \frac{\partial}{\partial x}f_{B}^{r,\mu} = \int_{
\mathbb{R}_{v_{\ast}}^{3} \times \mathbb{S}^{2}} \chi^{r}
\mathcal{B}_{\mu}^{BB}(v-v_{\ast},\omega)f_{B}^{r,\mu}(x,v')
f_B^{r,\mu}(x,v_{\ast}^{'})dv_{\ast} d\omega \qquad \qquad
\\ \nonumber + \int_{
\mathbb{R}_{v_{\ast}}^{3} \times \mathbb{S}^{2}} \chi^{r}
\mathcal{B}_{\mu}^{AB}(v-v_{\ast},\omega)f_{A}^{r,\mu}(x,v')
f_B^{r,\mu}(x,v_{\ast}^{'})dv_{\ast} d\omega \qquad \qquad
\\ \nonumber    - f_{B}^{r,\mu} \int_{ \mathbb{R}_{v_{\ast}}^{3} \times \mathbb{S}^{2}}\chi^{r} \mathcal{B}_{\mu}^{BB}(v-v_{\ast},\omega)
f_B^{r,\mu}(x,v_{\ast})dv_{\ast}d\omega
\\ \nonumber - f_{B}^{r,\mu} \int_{ \mathbb{R}_{v_{\ast}}^{3} \times
  \mathbb{S}^{2}}
\chi^{r} \mathcal{B}_{\mu}^{BA}(v-v_{\ast},\omega)
f_A^{r,\mu}(x,v_{\ast})dv_{\ast}d\omega, \hspace*{5 mm} (x,v) \in (-1,1)
\times \mathbb{R}_{v}^{3},
\\ \nonumber   f_{B}^{r,\mu}(-1,v) = \sigma(-1)   \lambda ' M_{-}(v) ,
\hspace*{1 mm}
\xi>0 ,
 \hspace*{4 mm} f_{B}^{r,\mu}(1,v) = \sigma(1) \lambda ' M_{+}(v) ,
 \hspace*{1 mm}
\xi<0 ,
\\
\end{eqnarray}
\noindent with
\begin{eqnarray*}
  \int min( \mu , (1 + | v | )^{ \beta}) f_{B}^{r,\mu}(x,v)dx
dv=1 .
\end{eqnarray*}
Here,
\begin{eqnarray*}
 \sigma(-1) = \frac{\int_{ \xi<0}| \xi |f_{B}^{r,\mu}(-1,v)dv }{ \int_{\xi>0}\xi f_{B}^{r,\mu}(1,v)dv+{ \int_{ \xi<0}| \xi |f_{B}^{r,\mu}(-1,v)dv}}
\end{eqnarray*}
and
\begin{eqnarray*}
 \sigma(1) = \frac{ \int_{ \xi>0} \xi f_{B}^{r,\mu}(1,v)dv }{ \int_{\xi>0}\xi f_{B}^{r,\mu}(1,v)dv+{ \int_{ \xi<0}| \xi |f_{B}^{r,\mu}(-1,v)dv}}.
\end{eqnarray*}
By using the mass conservation as in (\cite{B1}) we can prove that 
the boundary conditions of (\ref{rb}) writes 
\begin{eqnarray} \label{tr}
\nonumber f_{B}^{r,\mu}(-1,v) =  M_{-}(v) \int_{ \xi<0}| \xi
|f_{B}^{r,\mu}(-1,v)dv   , \hspace*{2 mm}
\xi>0,
\\  f_{B}^{r,\mu}(1,v) = M_{+}(v)  \int_{ \xi>0} \xi
f_{B}^{r,\mu}(1,v)dv , \hspace*{2 mm}
\xi<0.\qquad \qquad
 \end{eqnarray}
From here the arguments of (\cite{LA1}, \cite{LA2}, \cite{B1} \cite{Bw}) can
be used to pass to the limit in the parameters $(r , \mu) $ and to
prove 
that $(f_A,f_B)$ satisfies (\ref{cA}, \ref{cB}) in the weak sense for 
$0 \leq \beta <2$ and in the renormalized sense for $-3 < \beta \leq
0$.

But for the sake of clarity we explain the passage to the limit in the terms (\ref{tr}) i.e we prove the weak convergence in
 $L^{1}( \{ v \in \mathbb{
R}_{v}^{3} , \xi >0 \})$ ( resp $L^{1} ( \{ v \in
\mathbb{R}_{v}^{3} , \xi <0 \}) $) of $f_{B}^{j}(1,.)$ ( resp.
$f_{B}^{j}(-1,.)$) to $f_{B}(1,.)$ (resp. $f_{B}(-1,.)$). First, it is 
important to check that 
the fluxes $ \int_{\xi >0} \xi f_{B}^{j}(1,v)dv$ and $ \int_{\xi
<0} |\xi| f_{B}^{j}(-1,v)dv$ are controled. 
From (\ref{rb}) written in the exponential form, it holds that 
\begin{eqnarray} \label{ey}
\nonumber f_{B}^{j}(x,v)  \geq \hspace*{75 mm}
\\  \nonumber f_{B}^{j}(-1,v) e^{-
\int_{-\frac{1+x}{\xi}}^{0} \int_{\mathbb{R}_{v_{\ast}}^{3} \times
\mathbb{S}^2 }\chi^{r}(\mathcal{B}_{BA}^{\mu}f_A^{r,\mu}(x+s \xi,v_{\ast} )
+\mathcal{B}_{BB}^{\mu}f_B^{r,\mu}(x+s \xi,v_{\ast} )dv_{\ast}d\omega ds} ,
\\ \nonumber \xi >\frac{1}{2} , |v| \leq 2, \hspace*{5 mm}
\\ \nonumber f_{B}^{j}(x,v)  \geq  \hspace*{75 mm}
\\ \nonumber f_{B}^{j}(1,v) e^{-
  \int_{\frac{1-x}{\xi}}^{0} \int_{\mathbb{R}_{v_{\ast}}^{3} \times \mathbb{S}^2
  } \chi^{r}( \mathcal{B}_{BA}^{\mu}f_A^{j}(x+s \xi,v_{\ast})dv_{\ast} 
+ \mathcal{B}_{BB}^{\mu}f_B^{j}(x+s \xi,v_{\ast} )dv_{\ast} d\omega ds  }, 
\\ \xi < - \frac{1}{2} , |v| \leq
2.\hspace*{5 mm}
\end{eqnarray}
 For $v$ satisfying $ |v| \leq 2$ with $\xi > \frac{1}{2}$ or
$\xi < -\frac{1}{2}$, 
\begin{eqnarray*}
\int_{-1}^{1}   \int_{\mathbb{R}_{v_{\ast}}^{3} \times \mathbb{S}^2
  }   \frac{\chi^{r}}{|\xi|}   \left( \mathcal{B}_{BA}^{\mu} 
f_{A}^{r,\mu}(z,v) + \mathcal{B}_{BB}^{\mu} 
f_{B}^{r,\mu}(z,v) 
  \right)dv_{\ast} d\omega
dz
 \end{eqnarray*}
is uniformly bounded from above. Hence, using the definition
of the boundary conditions (\ref{bonb}) in (\ref{ey}),
\begin{eqnarray*}
f_{B}^{j}(x,v)& \geq & c M_{-}(v) \int_{\xi <0} |\xi|
f_{B}^{j}(-1,v)dv \quad , \quad \xi >\frac{1}{2} , \quad |v| \leq
2,
\\ f_{B}^{j}(x,v) & \geq & c M_{+}(v) \int_{\xi >0}\xi f_{B}^{j}(1,v)dv  \quad ,
\quad \xi < - \frac{1}{2} , \quad |v| \leq 2 .
\end{eqnarray*}
So,
\begin{eqnarray*}
 c \int_{\{\xi >\frac{1}{2} , |v| \leq 2 \} \cup \{\xi < -
\frac{1}{2}
  , |v| \leq 2 \} }f_{B}^{j}(x,v)dxdv
  \\  \qquad \qquad \geq \int_{\xi >0}\xi f_{B}^{j}(1,v)dv + \int_{\xi <0} |\xi| f_{B}^{j}(-1,v)dv .
 \end{eqnarray*}
 $f_{B}^{j}$ being non negative,
\begin{eqnarray*}
 c \int_{-1}^{1} \int_{ \mathbb{ R}_{v}^{3}}
min(\mu,(1+|v|)^{\beta})f_{B}^{j}(x,v)dx dv
\\  \qquad \qquad \geq \int_{\xi >0}
\xi f_{B}^{j}(1,v)dv + \int_{\xi <0} |\xi| f_{B}^{j}(-1,v)dv.
\end{eqnarray*}
Since $ \int_{-1}^{1} \int_{ \mathbb{ R}_{v}^{3}}
min(\mu,(1+|v|)^{\beta})f_{B}^{j}(x,v)dx dv =1 $, the fluxes
$\int_{\xi
>0} \xi f_{B}^{j}(1,v)dv$ and $\int_{\xi <0} |\xi| f_{B}^{j}(-1,v)dv$ are bounded uniformly
w.r.t $j$. 
\\Furthermore, the energy fluxes are also controlled. Indeed, from property \ref{cons}, the conservation of the
energy for $(f_A^{j},f_B^j)$, gives
\begin{eqnarray*}
 \int_{\xi >0}  \xi v^{2} f_{B}^{j}(1,v)dv + \int_{\xi <0}
|\xi|v^{2} f_{B}^{j}(-1,v)dv
 \\ \leq \int_{\xi >0}  \xi v^{2}
( m_A f_A^{j}(-1,v) + m_B f_B^j(-1,v) ) dv  
\\ + \int_{\xi <0} |\xi|v^{2} ( m_A f_A^{j}(1,v) + m_B f_A^{j}(1,v) ) dv.
\end{eqnarray*}
By definition of the boundary conditions (\ref{ra}) and
(\ref{rb}),
\begin{eqnarray} \label{ih}
  \nonumber  \int_{\xi >0}  \xi v^{2} f_{B}^{j}(1,v)dv + \int_{\xi <0}
|\xi|v^{2} f_{B}^{j}(-1,v)dv
\\    \leq  (k^{j} + \int_{\xi' <0} | \xi'|
f_{B}^{j}(-1,v')dv' )\int_{\xi > 0} \xi v^{2}M_{-}(v)dv
\\  \nonumber   +  (k^{j}
+\int_{\xi' >0}  \xi' f_{B}^{j}(1,v')dv' ) \int_{\xi <0} |\xi|v^{2}
 M_{+}(v)dv .
\end{eqnarray}
\noindent The right-hand side of (\ref{ih}) being bounded, the energy 
fluxes are also bounded. 
Finally, the entropy fluxes can also be controled. Indeed
\begin{eqnarray} \label{fb}
\nonumber \xi \frac{\partial}{\partial x} \big( f_A^{j} (log(f_A^{j}) -1)\big) =
Q_{AA}^j(f_A^{j},f_A^{j})log(f_A^{j}) + Q_{AB}^j(f_A^{j},f_B^{j})log(f_A^{j})
,
\\ \nonumber \xi \frac{\partial}{\partial x} \big( f_B^{j} (log(f_B^{j}) -1)\big) =
Q_{BB}^j(f_B^{j},f_B^{j}) \log(f_B^{j}) +
Q_{BA}^j(f_B^{j},f_A^{j})log(f_B^{j}).
\\
\end{eqnarray}
Using a Green's formula and an entropy estimate in the system (\ref{fb}), leads to
\begin{eqnarray*}
 \int_{\xi >0} \xi f_{B}^{j}(1,v) \log f_{B}^{j}(1,v)  dv +
\int_{\xi <0} |\xi| f_{B}^{j}(-1,v) \log f_{B}^{j}(-1,v)  dv
\\  \leq ( \int_{\xi ' >0} \xi' f_{B}^{j}(1,v ')dv' + k^{j})
\\  \int_{\xi <0} |\xi| M_{+}(v)
log(M_{+}(v)(\int_{\xi ' >0} \xi' f_{B}^{j}(1,v')dv' + k^{j})) dv
\\  + ( \int_{\xi' <0} |\xi'| f_{B}^{j}(-1,v')dv' + k^{j})
\\   \int_{\xi >0}  M_{-}(v) \log (M_{-}(v)(\int_{\xi' <0} |\xi'| f_{B}^{j}(-1,v')dv' + k^{j}))dv.
\end{eqnarray*}
\noindent By the Dunford-Pettis criterion
(\cite{CIL}), $f_{B}^{j}(1,.)$ is weakly compact in
\\ $L^{1}(\{ v \in
\mathbb{R}_{v}^{3}, \xi > 0 \}  )$. Let one of its subsequence still
denoted by $f_{B}^{j}(1,.)$, converging weakly to some $g_{+}$ in
$L^{1}( \{ v \in \mathbb{R}_{v}^{3}, \xi > 0 \}  )$. From now the
identification between $g_+$ and $f_B(1,v)$ is analogous to the proofs
given in (\cite{B1}, \cite{Bw}). This concludes the proof of Theorems
1 and 2.


\begin{thebibliography}{1}
\bibitem{A}
Aoki, K. The behaviour of a vapor-gas mixture in the
continuum limit: Asymptotic analysis based on the Boltzman equation. 
  T.J. Bartel, M.A.Gallis eds: AIP, Melville,
  2001, 565-574.
\bibitem{ABT}
Aoki K., Bardos C, Takata S. Knudsen layer for a gas mixture. 
Journal of Statistical Physics 2003, 112 (3/4); 629-655. 
\bibitem{ATT}
Aoki, K.; Takata, S.; Taguchi S. Vapor flows with evaporation
and condensation in the continuum limit: effect of a trace of non
condensable gas, 2003. European Journal of Mechanics B/Fluids 22, 51-71.
\bibitem{ATK}
Aoki, K.; Takata, S.; Kosuge, S. Vapor flows caused by evaporation
and condensation on two
  parallel plane surfaces: Effect of the presence of a noncondensable
  gas. Physics of fluids, 1998, volume 10, number 6, 1519-1532.
\bibitem{LA}
Arkeryd, L., Nouri, A., 'A compactness result related to the stationary Boltzmann equation in a slab, with applications to the existence theory', Ind. Univ. Math. Journ., 44 (3), 815-839, 1995.
\bibitem{LA1}
Arkeryd, L.; Nouri, A. The stationary Boltzmann Equation in
the Slab with Given
  Weighted Mass for Hard and Soft Forces. Ann.Scuola.Norm.Sup.Pisa 1998,
  27, 533-536.
\bibitem{LA2}
Arkeryd, L.; Nouri, A. $L^{1}$ solutions to the
stationary Boltzmann equation in a slab. Annales de la Faculte
des sciences de Toulouse 2000, 9, 375-413.
\bibitem{AM}
Arkeryd, L.; Maslova, N. On diffuse reflection at the
boundary for the Boltzmann equation and related equations. J.Stat.Phys. 1994, 77, 1051-1077. 
\bibitem{AH}
Arkeryd, L.; Heintz, A. On the solvability and
asymptotics of the Boltzmann
  equation in irregular domains. Comm.Part.Diff.Eqs. 1997, 22,
  2129-2152.
\bibitem{B1} 
S.Brull. \emph{The Boltzmann equation for a two component gas in the slab.} 
Math.Meth.Appl.Sci. 2008, 31, 153-178.
\bibitem{Bw} 
S.Brull. \emph{The Boltzmann equation for a two component gas in the
  slab for soft forces.} Math.Meth.Appl.Sci. (2008), 31, 1653-1666. 
\bibitem{B2} 
S.Brull. \emph{Problem of evaporation-condensation for a 
two component gas in the slab.} Kinetic and Related Models, (2008) Vol 1, No 2, 185-221.
\bibitem{CIL}
 Cercignani, C.; Illner, R.; Pulvirenti, M. Existence and Uniqueness Results. The mathematical theory of dilute gases. Springer:
  New York, 1994, 133-163.
\bibitem{chap} 
S.Chapman ; T.G. Cowling The mathematical theory of non uniform gases
(Cambridge univ. Press. , 1970). 
\bibitem{DPL} DiPerna, R.; Lions, P.L. On the Cauchy problem for 
the Boltzmann equation: Global existence and weak stability. 
Ann. Math. 130, 1989, 321-366. 
\bibitem{Do} C. Dogbe, Fluid dynamic limits for gas mixture. I. Formal derivations, Mathematical Models and Methods in Applied Sciences, Vol. 18,  no. 9, (2008), 1633-1672.
\bibitem{GT} Golse F., Takata, S. \emph{Half-space of the nonlinear 
Boltzmann equation for weak evaporation and condensation of a binary
mixture
of vapors}, Eur.J.Mech.B-Fluids 26 (1), 105-131 2007.
\bibitem{L}
Lions, PL. Compactness in Boltzmann equation via Fourier
integral
  operators and applications. J.Math.Kyoto Univ. 1994, 34,
  391-427.
\bibitem{M}
Mischler, S. On the weak-weak convergences and applications to the
initial boudary value problem for kinetic equations, preprint 2003.
\bibitem{P}
V. Panferov. On the interior boundary-value problem for the stationary Povzner equation with hard and soft interactions, Ann.Sc.Norm. Super. Pisa Cl. Sci. (5) 3 (2004), no. 4, 771-825.
\bibitem{SAD}
Sone, Y.; Aoki, K.; Doi, T. Kinetic theory analysis of gas flows
condensing on a plane condensed phase: Case of a mixture of a
vapor and noncondensable gas. Transport theory and statistical
physics, 1992, 21(4-6), 297-328.
\bibitem{SATSB}
Sone, Y.; Aoki K; Takata S; Sugimoto H.; Bobylev A Inappropriateness
of the heat conduction equation for description of a temperature field
off a sationary gas in the continuum limit: examination by asymptotic
and numerical computation of the Boltzmann equation Physics of fluids,
1996, volume 8, number 2, 628-638. Erratum: Physics of fluids 1996; 8,
841. 
\bibitem{T} 
Takata S.Kinetic theory analysis of the two-surface problem of
vapor-vapor mixture in the continuum limit, Physic of Fluids, 16, 7, 2004. 
\bibitem{TAL} 
S. Taguchi, K. Aoki, V. Latocha, Vapor flows along a plane condensed phase with weak condensation in the presence of a noncondensable gas, Journal of Statistical Physics, 124(2006) pp 321-369
\bibitem{TAT}
Taguchi, S.; Aoki, K.; Takata S. Vapor flows condensing at
incidence onto a plane condensed phase in the presence of a
noncondensable gas.I. Subsonic condensation. Physics of fluids, 2003, volume 15, number 3, 689-705.
\end{thebibliography}
\end{document}